\documentclass[11pt, leqno]{amsart}

\usepackage{amssymb}
\usepackage{amsmath}
\usepackage{amscd}
\usepackage{float}
\usepackage{graphicx}
\usepackage{amsfonts}
\usepackage{color}
\newtheorem{definition}{Definition}

\newtheorem{theorem}{Theorem}

\newtheorem{rem}{Remark}

\def\C{\mathcal C}

\def\x{\mathsf{x}}
\def\y{\mathsf{y}}
\def\K{\mathsf{K}}
\def\OO{\mathsf{O}}

\begin{document}

\title{AN INVARIANT OF  COLORED  LINKS VIA SKEIN  RELATION}

\author{Francesca Aicardi }
 \address{The  Abdus Salam  Centre for Theoretical Physics (ICTP),  Strada  Costiera, 11,   34151  Trieste, Italy.}
 \email{faicardi@ictp.it}

\dedicatory{To  the  memory  of  Serezha Duzhin}



\date{}

 \begin{abstract}
In  this note we  define  a polynomial  invariant  for  colored links by a  skein  relation.  It   specializes to  the   Jones  polynomial  for  classical  links.
\end{abstract}

 \maketitle
\section*{Introduction}
This  note  deals with an invariant of colored links.  Let us detail   these  objects.

 An  oriented  link   with  $n$  components is  a  set  of $n$  disjoint  smooth oriented  closed   curves embedded  in $S^3$.  An oriented  link  is  {\sl  colored},  if each  component  is  provided  with  a  color. In other  words,  a  link is  colored  if a  function $\gamma$ is  defined  on the  set  of  components  with values in a  finite   set  $N$  of  colors. Let's call   such a   function $\gamma$ a {\sl  coloration}.   Every  coloration introduces an  equivalence relation  in  the  set  of  components: two  components  are  equivalent  if  they  have  the  same  color.   Let $C$  be  the  set  of  components  of a  link $L$.  Two  colorations $\gamma$  and  $\gamma'$  of  $L$    are  said  to be {\sl  equivalent} if  there is  a  bijection in $N$  between  the   images   $\gamma(C)$  and  $\gamma'(C)$.  Two    colorations  define   the same  partition  of  the  set  $C$  into  classes  of equivalence if  and  only  if  they  are  equivalent.

An invariant of colored links takes the same value on isotopic links with the  same  coloration, and may take different values on the same link with different colorations.
Our  invariant of  colored links    takes the  same value on  isotopic  links  with    equivalent  colorations.

There  is a  wide literature  on  the  invariants  of  colored  links,  developed   from   the  multivariable  Alexander  polynomial, which includes   contributions  by Fox,  Conway,  Kauffman,   and  many others. For  a  summary  of  the history of these  invariants  see \cite{Cim1}. More  recent  works in  this  direction \cite{Mur}\cite{Cim1}\cite{Cim2} use also  skein  relations for  the  Alexander  polynomial,  but the  colors  of  the components  of  the links  involved  in  such  skein  relation  are  unaltered.

The skein  relation  introduced  in  the  present  note  may change  the  colors  of  the  components  of  the  link.  Therefore the  invariant  that  we  define  seems  to  be  new;  still, the  relation  of  it  with  the  known  invariants  has  to  be  investigated.

We recall  here  the  definition  of  an  invariant  of  classical  links  by  a  skein  relation.

Among  the most famous topological invariants  of oriented  knots  and  links, we single out   the  Alexander  polynomial, introduced  in  \cite{Ale}, and  the  Jones  polynomial, introduced in   \cite{joAM}. Both  these  polynomials in one  variable  can  be considered  as  particular  instances  of  a  more general  invariant  of  isotopy  classes  of  oriented  links,  the so-called   HOMFLY polynomial, introduced in \cite{Hom},  which  is  a Laurent  polynomial  in  two  variables.  This polynomial is  uniquely  defined  by  the  condition  that  its  value  is  1 on  the  unknotted  oriented  circle,  and  by a  skein  relation, as we  explain  below.

Suppose  that  $K_\x^+$,  $K_\x^-$  and  $K _\x^0$  are  diagrams  of  three  oriented  links  that are  exactly  the  same  except in  the  neighborhood  of a crossing point $\x$,  where  they look  as  shown in the  figure below.  Then  the skein  relation  is  an  equation  that  relates the  values  of  the  polynomial $P$ on  these  links:
\begin{equation}\label{skeinrule}   \ell P(K_\x^+)  +  \ell^{-1} P(K_\x^-) +  m  P(K_\x^0)= 0, \end{equation}

\centerline{\includegraphics{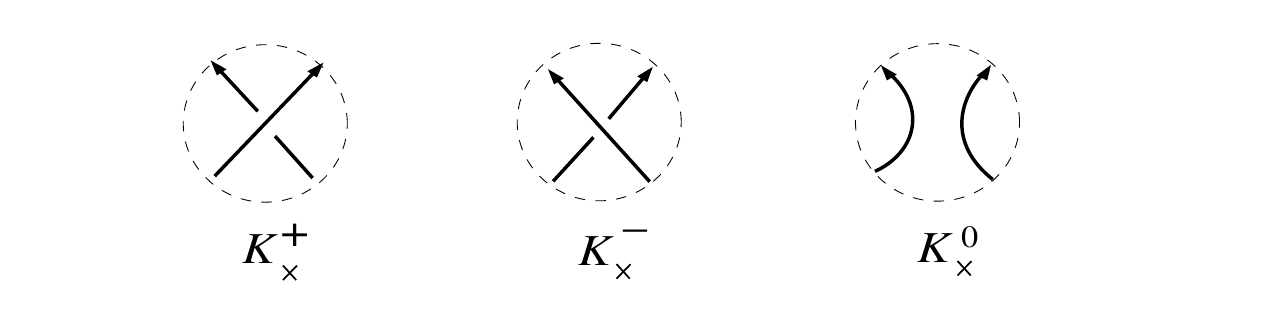}\label{color0}}

 The  proof that  the  skein  relation  allows  to define  the  value  of  $P$ on  any link,  is  done  by  induction  on  the  number  of  crossings.
By  an  ordered sequence    of  crossing  changes,  any  link $K$    can  be  changed into  a  trivial  link $K_u$,  i.e.  a  collection  of unknotted  and  unlinked  components.  The     value  of  $P$    on  $K_u$ is  supposed known: indeed,  it  must  be independent of the  number $n$  of  crossing points,   depending only  on  the  number of  components.  Let  $(\x_1,\dots, \x_m)$ be  the sequence of crossings  where the  sign has  to  be  changed to  transform  $K$ into  $K_u$,  and  let   $_iK$ be the  diagram  obtained  from  $K$  by  changing  the  sign  of  the  first $i$ crossings, so  that  $_mK=K_u$.  The  skein  relation  is thus  applied to the  crossing $\x_1$:  if  this  crossing  is  positive in  $K$,  then  $K=K_{\x_1}^+$, and  one  obtains $P(K_{\x_1}^+)$ as a linear  combination  of the values  of  $P(K_{\x_1}^-)$  and  $P(K_{\x_1}^0)$. Since  $K_{\x_1}^0$ has $n-1$ crossings  and  then  $P(K_{\x_1}^0)$ is  known  by  the induction  hypothesis.   Now, by  definition  $K_{\x_1}^-=\ _1K$, and  the value $P(_1K)$  is  calculated  by  using  the  skein  relation  applied  to the second crossing point  $\x_2$,  and  so  on, until the  last  crossing  change, in  which  either $_{(m-1)}K_{\x_m}^+$  or  $_{(m-1)}K_{\x_m}^-$ coincides with  $K_u$.  Then     the calculation  of     $P(K)$  is  concluded,  being    obtained  from   $P(K_u)$  and  the values of  $P$ on $(m-1)$  links  with  $n-1$  crossings.

 It is  thus  proved  that $P(K)$  calculated   in this  way is  independent  of  the  chosen  sequence  of  changes,  as  well as independent of  the  particular  projection  of  the  link.   The  polynomial  $P$   is  in  fact  an  isotopy  invariant  of  the  link.

Suppose  now  that  we  want to  define  by  a  skein  relation  an invariant, which is able  to distinguish  two  links  that  are  identical,  but have  non-equivalent  colorations.
We  show  that  it  is possible  to  obtain the  invariant  we  want,  by introducing   a  new  skein  relation,  which  takes  into  account   strands  of  any  color,  and by setting  the  value of  the  invariant  on  the  trivial  links with  non-equivalent    colorations,  i.e.,   on  collections  of $n$  unknotted  and  unlinked  components  that  are   colored,  for  every $c=1,\dots,n$,   with  $c$  different  colors.

Our  invariant, $F$, requires  three  variables,  say  $x$,  $w$  and  $t$.  When  the link  is  colored with  a  sole  color,  our invariant  is     reduced  to another  instance of the  HOMFLY   polynomial,  and  becomes  the Jones  polynomial  in  the  variable  $t$  if  $w=t^{1/2}$.

Finally, we remark   that  the  polynomial  $F$  was  not  found by searching for  an  invariant  for  colored  links.  Its definition  is  a  byproduct  of  the  study of another  class  of  links,  the  {\sl  tied  links} (see \cite{AJTL}),  obtained  from  the  closure of  tied  braids,  which  in  turn  constitute  a  diagrammatic  representation  of  an  abstract  algebra,  the  so-called  nowadays   {\sl algebra  of  braids  and  ties}. The defining set of generators of the algebra of braids and ties consist of two sets of generators, one of them is interpreted as usual  braid  generators and the other one as the  set of ties. Despite  the fact that the defining  relations involving  ties generators     look complicated and,  in  particular, the ties do  not  commute  with  the  braid elements, we have  shown \cite{AJtrace}  that a new   geometrical  interpretation  of the ties-generators,   coherent  with  the defining relations of the algebra,   allows  in particular  to reduce  the  ties  to  simple  connections  between different  strands  of  the  braids,  and,  under  closure,  between  different  components  of  the  links.  These  connections  behave  as no other  than  an  equivalence  relation  between  link  components. As  a  consequence,  the invariant  polynomial for  tied  links  defined  in    \cite{AJTL} provides an  invariant  of  colored  links  as  well.

\section{An invariant for  colored  links}\label{Sec2}

A  colored link diagram  looks like a  diagram of a  link,  where  each component  is  colored.

\begin{definition}\label{isotop} \rm Two oriented colored links    are  {\sl c-isotopic}  if
  they  are  ambient  isotopic,  and  their  colorations  are  equivalent.

\end{definition}

\begin{figure}[H]
\centering \includegraphics{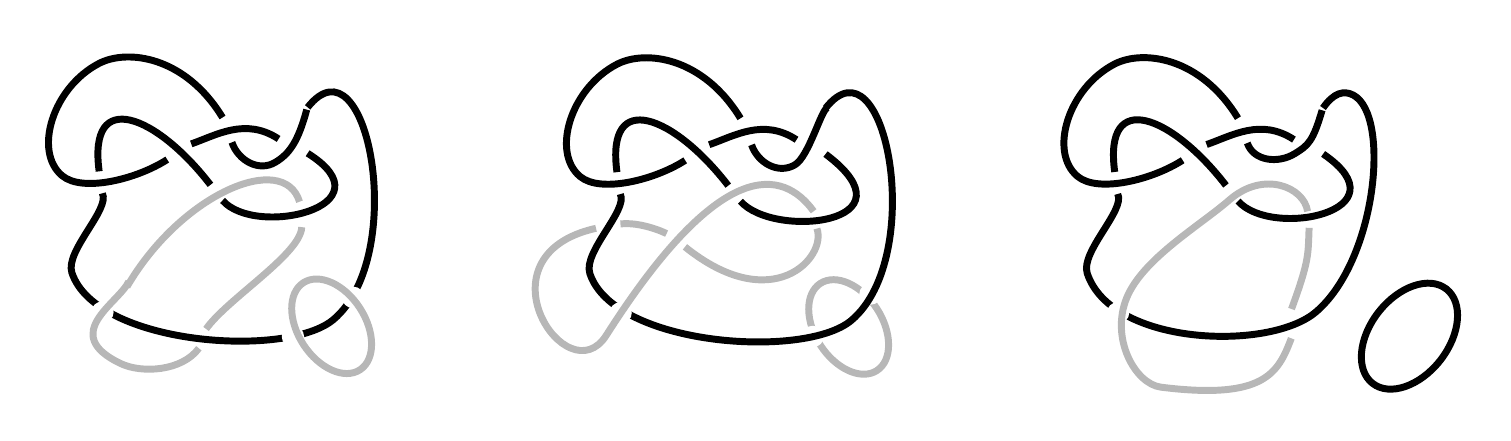}
\caption{The  first  two  diagrams represent  c-isotopic colored links,  the third  one  is  not c-isotopic  to  the  others. Black and  gray  represent     two different colors.}
\end{figure}


Let  $R$  be  a  commutative  ring and  $\C$  be  the set of   colored oriented link  diagrams.
By  an {\sl invariant of  colored  links}, we  mean     a  function $I:  \C \rightarrow R$, which is  constant      on each class  of c-isotopic links.

 \begin{rem}\rm  In  the  sequel, we denote  by  $\K$    an  oriented   colored  link  as  well as  its  diagram,  if  there  is  no risk  of  confusion. \end{rem}

\begin{figure}[H]
\centering
\includegraphics{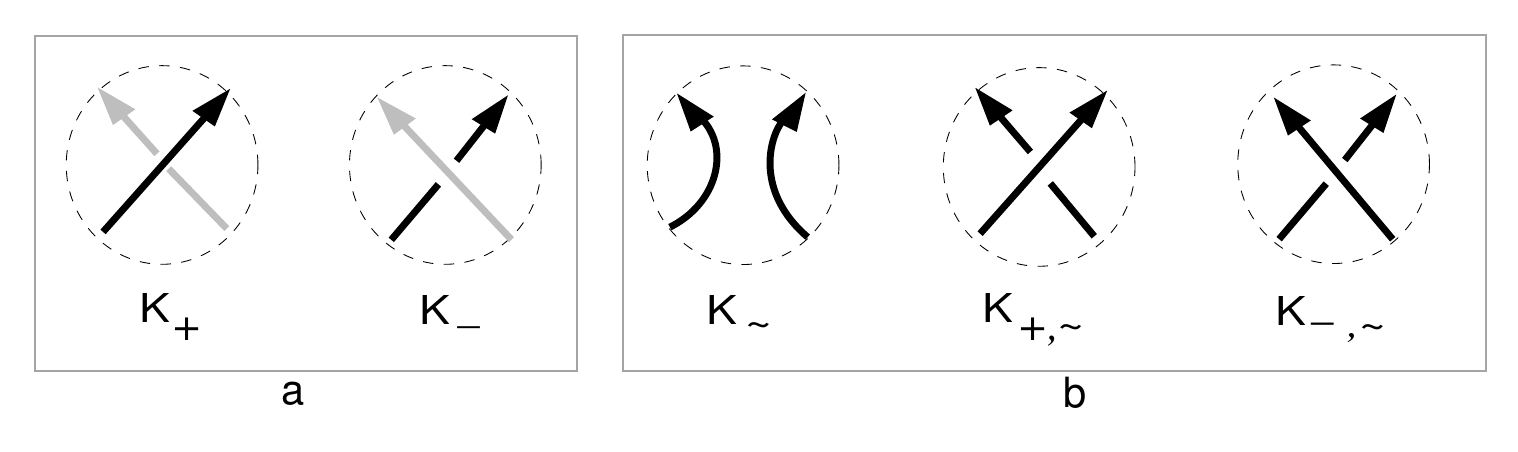}
\caption{a) The  discs,  where   $\K_+$  and $\K_-$  are  not c-isotopic. b) The  discs where  $\K_{\sim}$,   $\K_{+,\sim}$,  and $\K_{-,\sim}$   are not c-isotopic.  Black  and  gray  indicate  any  two   colors. }\label{color2}
\end{figure}

 \begin{theorem} \label{the1} There  exists  a rational function  in  the  variables $x,w,t$,   $F:\C \rightarrow   \mathbb{Q}(x,t,w)$,     invariant  of oriented colored  links,    uniquely  defined  by  the   following  three  conditions  on  colored-link  diagrams:
\begin{itemize}
\item[I] The  value of  $F$ is  equal to 1  on  the  unknotted circle.

\item[II] Let  $\K$   be  a colored link,  with  $n$  components  and  $c$  colors.  By  $[\K;\OO]$  we  denote  the  colored link  with  $n+1$ components consisting of $\K$  and  the unknotted and unlinked  circle, colored  with  the ($c+1$)-th     color.
Then

\[  F([\K;\OO])=  \frac {1}{wx} F(\K). \]

\item[III](skein  relation)  Let $\K_+$  and   $\K_-$   be  the  diagrams  of two colored  links,  that are  the  same outside    a  small  disc into  which two  strands  enter,   and  inside  this  disc  look  as
   shown  in figure \ref{color2}a. Black and gray  indicate  any two  colors,  not  necessarily  distinct.  Similarly,  let   $\K_\sim$  and   $\K_{+,\sim}$  be two links  that  inside  the  disc look  as  shown  in  figure \ref{color2}b,  while outside  the  disc  coincide with $\K_+$  and  $\K_-$, except, possibly,  for  the  colors of all the components of the  link having  the same colors as    the  strands entering  the  disc:  all these  components   in   $\K_\sim$  and   $\K_{+,\sim}$  are  colored by  a  sole  color. Then  the  following  identity  holds:
 \[ \frac{1}{  w}F( \K_ +)-w F (\K_-) = \left(1-t^{-1}  \right) F (\K_\sim) + \frac{1}{w} (1-t^{-1} ) F(\K_{+,\sim}).\]

\end{itemize}
\end{theorem}

\begin{rem}\label{rem2}\rm

The skein relation  III  holds  for  {\sl any}  two  colors  of  the      strands involved. In  particular, if  the  two  colors coincide, then  $\K_+=\K_{+,\sim}$  and  $\K_-=\K_{-,\sim}$, so  that relation III  is  reduced  to  the relation

 \begin{itemize}
 \item[IV]
 \[ \frac{1}{t w }F( \K_ {+,\sim})- w F (\K_{-,\sim}) = (1-t^{-1}) F (\K_\sim).\]
\end{itemize}
Relation III also implies  the following  two  relations; they  will  be  used in the  sequel.

\begin{itemize}
\item[Va]
 \[ \frac{1}{ w }F( \K_ {+})= w \left[F(\K_-) +\left( {t}-1 \right) F(\K_{-,\sim})\right] + \left(t-1\right) F (\K_\sim).\]

\item[Vb]
\[   w  F( \K_ {-})= \frac{1}{ w } [F(\K_+) +(t^{-1}-1) F(\K_{+,\sim})] + ({t^{-1}}-1) F (\K_\sim).\]
 \end{itemize}

\end{rem}

\begin{rem}\label{rem3} \rm As  an  example of using  relations I, II  and  III,  let's  calculate  the  value  of  $F$ on  the  unlink $\OO^c_n$, consisting of  $n$ circles  and  $c$ colors.
If  $n=c$, then,  by II, $$F(\OO^c_c)= \frac{1}{wx}F(\OO^{c-1}_{c-1})=\left(\frac{1}{wx} \right)^{c-1}.$$
If  $n>c$,   then there are two circles colored by the same color. Regard them as $\K_{_{\displaystyle\widetilde{\ \ }}}$. Then $\K_+=\K_-=\OO_{n-1}^c$ (see  figure).

\centerline{\includegraphics{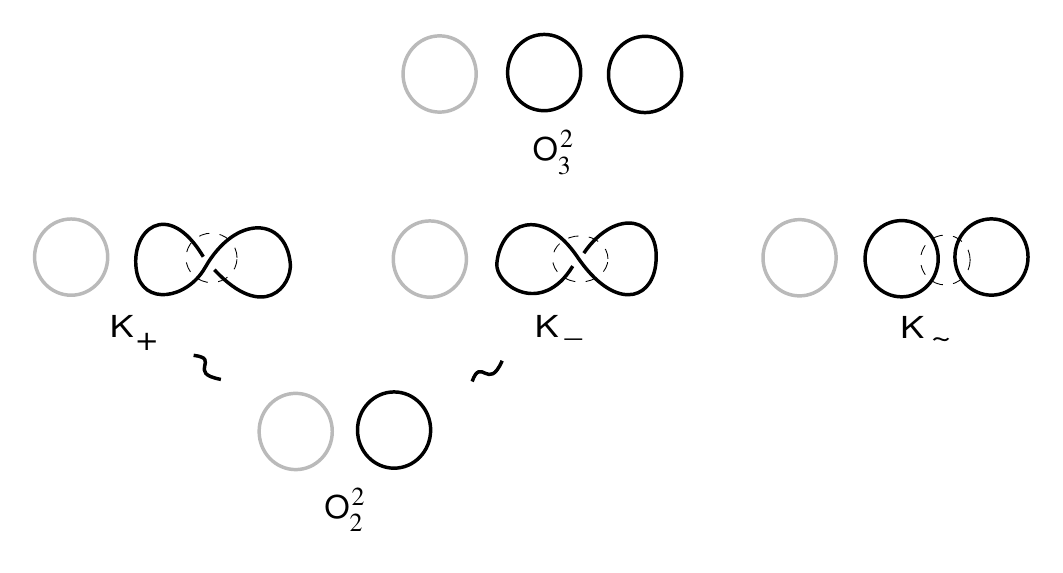}\label{color0}}

By IV,
$\displaystyle (1-t^{-1})F(\OO_n^c)=\Bigl(\frac1{tw}-w\Bigr)F(\OO_{n-1}^c)$.
So
$$F(\OO_n^c)=\frac{1-tw^2}{tw(1-t^{-1})}F(\OO_{n-1}^c) =
\frac{tw^2-1}{(1-t)w}F(\OO_{n-1}^c)= \frac{y}{xw}F(\OO_{n-1}^c),
$$
 where  \begin{equation}\label{formulay}  y=  \frac {x (t w^2-1)}{1-t}.\end{equation}
  By induction,
$\displaystyle F(\OO_n^c)=\Bigl(\frac{y}{xw}\Bigr)^{n-c}F(\OO_c^c)$.
 Thus
 \begin{equation}\label{formulaonc} F(\OO_n^c)=\frac{y^{n-c}}{(wx)^{n-1}} . \end{equation}

\end{rem}

{\bf Proof of Theorem 1.}
Theorem 1 is  proved by  the  same  procedure used in  the proof of the corresponding  theorem for classical links stated  at  page  112  in  \cite{Lic}.  We  will outline  the  parts  where  the  presence  of  colors modifies  the  demonstration.

Of  course,  the  Reidemeister  moves for  colored  links  are  all  the Reidemeister  moves,  where  the  strands  involved  may  have  different  colors.

  One starts with zero crossings as the induction base: the value of $F$ on the unlink $O_n^c$ with $n$ circles and $c$ colors is given by
$\frac{y^{n-c}}{(wx)^{n-1}}$ in accordance with Remark \ref{rem3}.

 Let  $\C^p$  be the  set  of  diagrams  of oriented colored  links  with $p$  crossings,  and  $\K\in \C^p$. By ordering  the  components  and  fixing a  point     in  each  component, one  constructs, for  every  diagram   $\K$,  an  associated  standard ascending diagram   $\K'$ in  the following   way:  traverse the  components  of  $\K$ in  their  given  order  and  from their  base points  in  the  direction  specified  by  their  orientation.   Every    crossing, encountered the first  time,  is  either  over or  under crossing. In  the  first  case, the  crossing   is  changed; otherwise, it  is  left as  in $\K$.  $\K'$  consists  of  the  same  number  of  components as  $\K$,  completely  unknotted  and  unlinked.  $\K$ and   $\K'$  are  identical,  except    for a  finite  sequence  of  crossings,  which we call 'deciding',  where  the signs  are  opposite.    Furthermore,  we  define    $\tilde \K'$,  as obtained  from  $\K'$ by  coloring with  a sole  color  all the  components having the same two  colors as  the two strands involved  into  {\it each}  deciding crossing.   $\K'$   and  $\tilde \K'$ are,  by  construction,    collections  of  unknotted  and  unlinked  components;  $\K'$ has   the  same number  of  colors   as  $\K$,  $\tilde \K'$  may  have a smaller number of  colors.  The  procedure  defining  $\K'$ allows  to  get  an  ordered   sequence  of deciding  crossings, whose  order  depends  on  the ordering  of  the
   components,    and  on the  choice of the base points.

 The  induction  hypothesis   states  that,  for each oriented  colored link  diagram $\K$ in  $\C^p$  there  is associated  a function  $F(\K)$, satisfying  relations I -- III,  which  is  independent of  the  ordering    of  the  components,   independent  of  the  choices  of  the  base  points, and invariant  under  Reidemeister  moves  that  do  not  increase  the  number  of  crossings  beyond  $p$.  Moreover, the  induction hypothesis states  that the  value  of  $F$  on any   colored link with  $p$  crossings,
  consisting   of  $n$  components  unknotted and unlinked,  and colored  with    $c$  colors  ($c\le n$),
  is the  same  as  $F(\OO_n^c)$ given  by  \eqref{formulaonc},   independent  of $p$.

 Now, let  $\K$  be  in  $\C^{p+1}$.   If  $\K$ consists   of  $n$  unknotted  and unlinked components, colored with  $c$ ($c\le n$) colors,  let's define
 \begin{equation}\label{eq2} F(\K)=   \frac{y^{n-c}}{(wx)^{n-1}}.  \end{equation}   Otherwise,  construct  the associate  ascending  diagrams $\K'$ and  $\tilde \K'$,  and
 consider  the  first    deciding  crossing $\x$.   If  in  a  neighborhood of  $\x$  the  colored  link  looks  like $\K_{+,\sim}$ (or  $\K_{-,\sim}$),
  use  skein  relation  IV    to  write   the  value  of  $F$  in  terms  of     $\K_{-,\sim}$ and  $\K_{\sim}$ (respectively, $\K_{+,\sim}$ and $\K_{\sim}$). If  in  a  neighborhood of  $\x$  the  colored  link  looks  like $\K_+$  (respectively, $\K_-$)   use  skein  rule Va  (respectively,  Vb)  to  write    the  value  of  $F$ in  terms of  the    value of $F$  on  the  colored  links    $  \K_{-},  \K_{-,\sim}$ and  $\K_{\sim}$ (respectively, $  \K_{+},  \K_{+,\sim}$ and  $\K_{\sim}$).  Observe  that  if  the  diagrams $\K_{\epsilon}$  or  $\K_{\epsilon, \sim}$ ($\epsilon=\pm$)   coincides,   in a  neighborhood  of  the  crossing  $\x$,    with  the  original  diagram, then $\K_{-\epsilon}$ or  $\K_{-\epsilon,\sim}$  coincides, in the  neighborhood of  the same  crossing,   with  the associated  diagrams  $\K'$ or $\tilde \K'$,  respectively.  On  the  other  hand, $\K_{\sim}$ represents  a  colored link diagram  with  $p$  crossings,  for  which  the  value  of  $F$  is  known,  and  invariant  according  to  the  induction  hypothesis. Then  we  look  for
   the  second  deciding crossing,  which  is  present  in  all the  diagrams  having  $p+1$ crossings,  obtained  by applying   the  skein  relation to  $\K$  at  $\x$. We  apply the  same  procedure to such diagrams  at  the  second deciding  crossing,  and  so  on.   The  procedure  ends  with  the  last  deciding  crossing, thus yielding  diagrams of    unlinked  colored  links, with  $p+1$ crossings,  where  the  function is given  by \eqref{eq2}, which depend  only  on  the  number  of  components  and  the  number  of  colors. For  an  example  see  subsection \ref{exa}.

  It  remains  to  prove  that

  \begin{enumerate}
  \item  the  procedure  is  independent  of  the  order of  the  deciding points;
  \item  the  procedure  is independent  of  the order of  components,   and of the choice  of  base-points;
  \item the  polynomial $F$
     is  invariant  under  Reidemeister  moves.

  \end{enumerate}

Following  the  proof done  in \cite{Lic} for  classical  links,  we  observe  that the  proofs  of statements (1),  (2)  and  (3) can be  done  in an  analogous way   also  in     presence  of colors.  Of  course,   every  time  a  skein  relation  is  used,  we  have  to  pay  attention   to  the  colors  of  all links  involved in.  We  write  here  the  proof  of  statement (1)  as  an  example.  The  proofs  of the other statements  are  similar.

The  proof  of statement (1)  consists  in  a  verification  that  the  value  of  the  invariant  does  not  change  if  we   interchange  any  two  deciding  crossings  in  the  procedure of  calculation.  So,  let  $\K$  be the  diagram  of  a  colored  link      and let  $\x$  and  $\y$  the  first two  deciding crossings  that  will  be  interchanged.

Denote by  $\epsilon_\x $  the   sign  at  the   crossing  $\x$, by $ \sigma_\x \K$ the  colored  link diagram  obtained  from $\K$  by changing  the  sign  at the  crossing   $\x$,  by $\tilde  \sigma_\x \K$ the diagram obtained from  $\K$ by changing  the  sign  at    $\x$ and coloring  by  a  sole  color all  the  components  colored by  the  two  colors of  the two  components crossing  at  $\x$;  then denote    by  $\rho_\x \K $ the  diagram  obtained  from $\K$  by removing  the  crossing     $\x$ and  coloring  by  a  sole  color all  the  components colored by  the  two  colors of   the  components  crossing  at $\x$.   Then do the same for the crossing   $\y$.

If    $\y$  follows  $\x$, then, by  skein  Va,b   :
\[ F(\K)=  w^{2\epsilon_\x}[ F ( \sigma_\x \K) +  (t^{\epsilon_\x}-1) F(\tilde  \sigma_\x \K)] + w^{\epsilon_\x}(t^{\epsilon_\x}-1 ) F(\rho_\x \K), \]
and
\[ F(\K)=  w^{2\epsilon_\x} \{ w^{2\epsilon_\y} [F( \sigma_\y \sigma_\x \K )    +  (t^{\epsilon_\y}-1) F(\tilde \sigma_\y \sigma_\x \K)]+w^{\epsilon_\y}(t^{\epsilon_\y}-1) F(\rho_\y\sigma_\x \K ) \}+ \]
\[  +(t^{\epsilon_\x}-1)w^{2\epsilon_\x}  \{ w^{2\epsilon_\y} [F( \sigma_\y \tilde \sigma_\x \K )    +  (t^{ \epsilon_\y}-1) F(\tilde \sigma_\y \tilde \sigma_\x \K)]+w^{\epsilon_\y} ( t^{ \epsilon_\y}-1 )F(\rho_\y\tilde\sigma_\x \K ) \}+ \]
 \[  +(t^{ \epsilon_\y}-1) w^{\epsilon_\x} \{ w^{2\epsilon_\y}[ F ( \sigma_\y \rho_\x \K) +  (t^{ \epsilon_\x}-1) F(\tilde  \sigma_\y \rho_\x \K)] + w^{\epsilon_\y}( t^{ \epsilon_\y}-1) F(\rho_\y  \rho_\x \K).\}  \]

If  $\x$  follows  $\y$,  then $F(\K)$ is  obtained  from  the  above  expression  by interchanging  $\x$ with  $\y$. Observe that this expression contains  terms  of type $\tau_\y \tau_\x \K$  or $\alpha(\tau_\x \tau'_\y \K +\tau'_\x \tau_\y \K)$,  where $(\tau,\tau')\in \{ \sigma, \tilde \sigma, \rho \}$ and  $\alpha$ is  a coefficient.  Such  terms  are   invariant  under  the  interchange of $\x$ with $\y$,  because  the  operation  $\tau_\x$  commutes with  $\tau_\y$  as  well  as with  $\tau'_\y$ (observe  that  the procedures  of coloring by  a  sole  color the  components  of the  two  colors  of  the  strands  crossings respectively  at $\x$ and  at $\y$, produce  equivalent colorations of  the  link   if they  are   interchanged).  Therefore  $F(\K)$ is  independent  of  the  order  of  $(\x,\y)$.

\hfill  $\square$

\subsection{Properties  of  the  Polynomial  $F$}

Here  we  list  some  properties  of  the  polynomial  $F$,  which  can  be  easily  verified.

\begin{enumerate}
\item[(i)]  $F$  is  multiplicative  with  respect to the   connected  sum of  colored  links.
\item[(ii)] The  value  of  $F$  does  not  change  if  the  orientations  of  all  curves  of  the  link  are  reversed.
\item[(iii)] Let $K$ be a  link  diagram whose components have all the same color, and $K^\pm$ be the  link  diagram  obtained  from  $K$  by changing  the  signs of  all  crossings. Thus   $F(\K^\pm)$  is  obtained  from $F(\K)$ by  the  following  changes:   $ w \rightarrow 1/w$  and    $t \rightarrow  1/t $.
\item[(iv)] Let  $\K$  be  a  knot or   a  link whose components  have  all the  same  color. Then    $F(\K)$  is  defined  by relations I and IV,  and   is  a  HOMFLY polynomial  with
     \[   \ell = i  \frac{1}{ w\sqrt{t}}   \quad \text{and} \quad m=i\left( \frac {1}{\sqrt t}-  \sqrt{t} \right). \]
     In  particular,  for  $w=t^{1/2}$  it  becomes  the  Jones  polynomial.
\end{enumerate}

Item (i)  follows  from  the  defining  relation  I of  $F$,  and  by  the  same  arguments proving  the  multiplicativity  of  the  invariants  obtained  by  skein  relations  (see  \cite{Lic}).  Items  (ii)  is  evident,  since the  value of  $F$  on  the  unlinked  circles is  independent  of  their  orientation,  and    the  skein  relations  are  invariant under  the  inversion  of  the  orientations  of  the  strands.  Item  (iii) follows from  the  fact that, if  $ w $ and  $t$ are  replaced  by $1/w$  and    $1/t $  respectively,  the  skein  relations   Va and  Vb are  interchanged,  whereas the term $y/(xw)$ stays unchanged.

As  for  item  (iv),  observe  that    comparing the   skein  relation  IV,  multiplied  by  $ t^{1/2}$,  with  equation (\ref{skeinrule}),     we  obtain the  expressions of  $\ell$  and  $m$  in  terms  of  $x,y$  and  $t$. Furthermore,  since for  the  Jones  polynomial $V(\K)$
\[ \ell=  i  t^{-1 }    \quad  \text{and}  \quad    m=i( t^{-1/2}-t^{1/2}), \]
 we  observe  that  the  expressions  of  $m$  are the  same  for $F(\K)$  and  $V(\K)$,  whereas   the  expressions  of  $\ell$  coincide  under  the  equality $w=t^{1/2}$.


\subsection{An example}\label{exa}

\begin{figure}[H]
\centering
\includegraphics{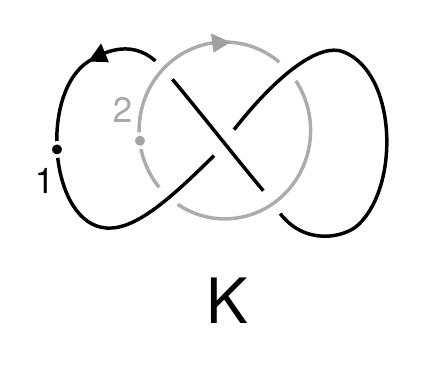}
\caption{The  components of  $\K$  are numbered  and  a base point is  chosen  on  each of them}\label{CL4}
\end{figure}

Let us apply  the  construction described  in  Theorem 1 to  calculate  the  polynomial  $F$ for the oriented colore link   $\K$    shown  in  figure \ref{CL4}.

\begin{figure}[H]
\centering
\includegraphics{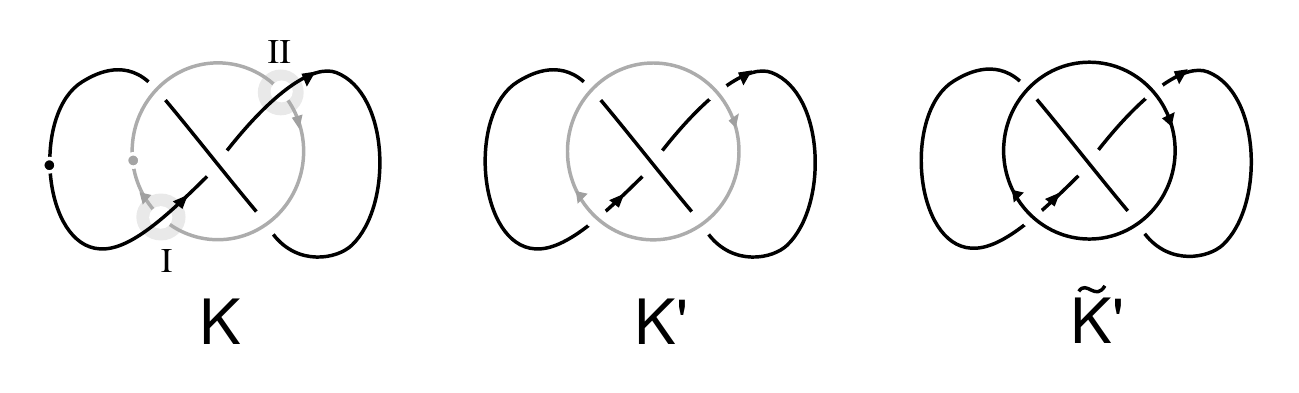}
\caption{The  deciding  crossings  are  marked as  I  and  II}\label{CL5}
\end{figure}

In  Figure, \ref{CL5} we  show  $\K$  with  the deciding  crossings  and  the  associated  diagrams  $\K'$  and $\tilde \K'$.

\begin{figure}
\centering
\includegraphics{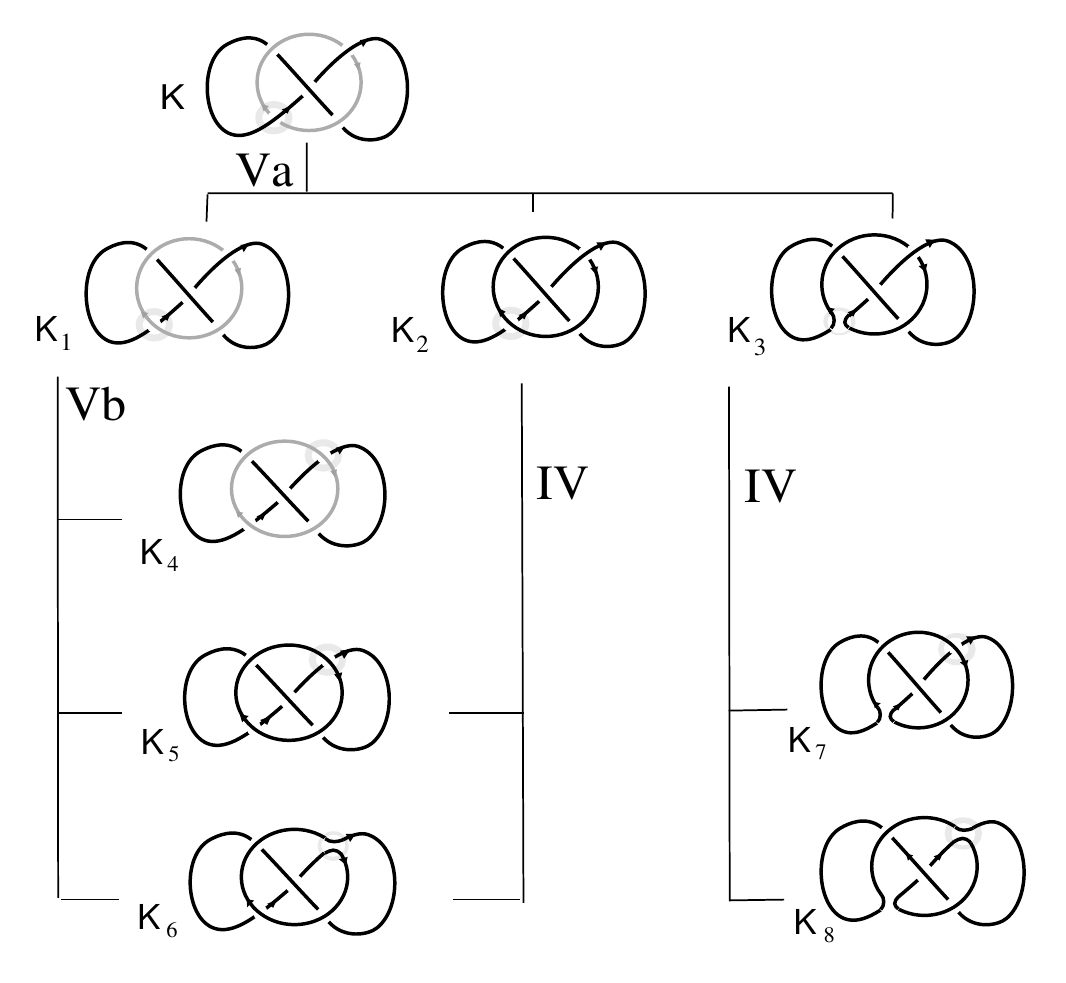}
\caption{}\label{CL6}
\end{figure}

The  first  deciding  crossing  is positive, the  second  one  is negative.  Applying  relation Va to  the  crossing I,  we  obtain (see  Figure \ref{CL6}):
\begin{equation}\label{K1}
F(\K)=w^2[F(\K_1)+ (t-1) F(\K_2)]+w(t-1)F(\K_3).
\end{equation}
Now we calculate    $F(\K_1)$ by  applying relation Vb  to  the  crossing II of $\K_1$:
\begin{equation}\label{K2}
 F(\K_1) =w^{-2}[F(\K_4) + (t^{-1}-1) F( \K_5)]+w^{-1}(t^{-1}-1)F(\K_6);
\end{equation}
then we  calculate 
$F(\K_2)$ by  applying relation IV  to  the  crossing  II of $\K_2$:
 \begin{equation}\label{K3}
 F(\K_2) =t^{-1}w^{-2} F(K_5)+  w^{-1}(t^{-1}-1)F(\K_6),
\end{equation}
  and  calculate $F(\K_3)$ by applying relation IV  to  the  crossing II  of $\K_3$ :
\begin{equation}\label{K4}
 F(\K_3) =t^{-1}w^{-2} F(\K_7)+w^{-1}(t^{-1}-1)F( \K_8).
 \end{equation}
Then  we  observe  that:

 $  \K_4 = \K'$,
and therefore, $F(\K_4)=1/(wx)$;

 $ \K_5=  \tilde \K'$, and therefore,   $F(\K_5)=y/(wx)$,  $y$  being   given  by (\ref{formulay});

$\K_6$  and   $\K_7$  are  both  c-isotopic to   $\OO$,
 the  unknotted  circle;  therefore    $F(\K_6)=F(\K_7)=1$;

 $ \K_8$  and $\K_2$  are c-isotopic; indeed,  they  are   different  projections  of  two oriented  circles of  the  same  color linked  by two  negative  crossings.  Therefore by  equation (\ref{K3}):

$$ F(\K_2)=F(\K_8)=   \frac{ w^2 t^2-  w^2 t+w^2-1  }{ w^3 t(1-t)  }. $$

Substituting  the  values  of  $F(K_7)$  and  $F(K_8)$ in equation (\ref{K4})  we  get

$$F(\K_3)=\frac{w^2 t^2+w^2-1  }{ w^4t^2}.$$

Since,  by  equation (\ref{K2}),  we  have

$$ F(\K_1)= \frac{w^2x+t-x }{w^3xt},$$

we  get, by equation (\ref{K1}):
\begin{equation}\label{FK}  F(\K)=  \frac{w^4(xt^2-xt^3)+w^2(t^2+xt-x-xt^2+xt^3)+x-xt}{w^3xt^2}.\end{equation}

Observe  also  that  $\K_3$ is c-isotopic  to  the oriented  trefoil  with  three  negative  crossings, and that  $\K_1$  is  c-isotopic  to  two oriented circles of  different  colors linked  by two  negative  crossings.


\begin{thebibliography}{25}
\bibitem{Cim1} D. Cimasoni, {\em A geometric construction of the Conway potential function}, Comment. Math. Helv. {\bf 79} (2004) 124–-146.

\bibitem{Mur} J. Murakami, {\em On local relations to determine the multi-variable Alexander polynomial of colored links},   Knots 90 (Osaka, 1990), de Gruyter, Berlin, 1992, 455–-464.

\bibitem{Cim2} D. Cimasoni and V. Florens, {\em Generalized Seifert surfaces and signatures of colored links}, Trans. Amer. Math. Soc.
{\bf 360}(3) (2008) 1223–-1264.

\bibitem{Ale}  J.W. Alexander, {\em Topological  invariants  of knots  and  links}, Trans.Am. Math. Soc. {\bf 30} (1928), 275--306.


\bibitem{joAM} V.F.R. Jones, {\it Hecke algebra representations of braid groups and link polynomials},
Ann. Math. {\bf 126} (1987), 335--388.



\bibitem{Hom} P. Freyd, D. Yetter, J. Hoste, W.B.R. Lickorish,  K.C. Millet, A. Ocneanu   {\em A new  Polynomial  Invariant  of  Knots  and Links.}  Bull. A.M.S. {\bf 12} (1985), 239--246.


\bibitem{AJtrace} F. Aicardi, J. Juyumaya, {\em Markov trace on the algebra of braids and ties.} submitted  to  Moscow Mathematical  Journal, october 2014, see  also arXiv:1408.5672.

 \bibitem{AJTL} F. Aicardi, J. Juyumaya, {\em Tied  Links.} Preprint ICTP,  november 2014, see  also arXiv:1503.00527.



\bibitem{Lic} W.B.R. Lickorish,  K.C. Millet,   {\em A Polynomial  Invariant  of  Oriented Links.}  Topology  Vol. 26,  No 1 (1987), 107--141.




\end{thebibliography}
\end{document}